\documentclass[final]{article}

\usepackage{amsmath}
\usepackage{amsfonts}
\usepackage{graphicx}
\usepackage{setspace}
\usepackage{citesort}
\usepackage{dcolumn}
\usepackage[squaren]{SIunits}
\usepackage{url}

\newcommand{\E}{\mathrm{e}}

\newcommand{\ba}{\boldsymbol{\alpha}}
\newcommand{\bx}{\mathbf{x}}

\newcommand{\field}[1]{\mathbb{#1}}

\newcolumntype{d}[1]{D{.}{\cdot}{#1}}

\begin{document}

\title{Moving least squares via orthogonal polynomials\thanks{Part of
    this work was carried out at the Institut f\"{u}r
    Str\"{o}mungsmechanik und Hydraulische Str\"{o}mungsmaschine of
    the University of Stuttgart, Germany, as part of the European
    Community-funded project HPC-Europa, contract number 506079}}

\author{Michael Carley\thanks{Department of Mechanical
    Engineering, University of Bath, Bath BA2 7AY, United Kingdom
    ({\tt{m.j.carley@bath.ac.uk})}}}

\bibliographystyle{siam}

\maketitle

\begin{abstract}
  A method for moving least squares interpolation and differentiation
  is presented in the framework of orthogonal polynomials on discrete
  points. This yields a robust and efficient method which can avoid
  singularities and breakdowns in the moving least squares method
  caused by particular configurations of nodes in the system. The
  method is tested by applying it to the estimation of first and
  second derivatives of test functions on random point distributions
  in two and three dimensions and by examining in detail the
  evaluation of second derivatives on one selected configuration. The
  accuracy and convergence of the method are examined with respect to
  length scale (point separation) and the number of points used. The
  method is found to be robust, accurate and convergent.
\end{abstract}

% \begin{keywords} 
%   moving least squares, interpolation, numerical differentiation,
%   orthogonal polynomials
% \end{keywords}

% \begin{AMS}
%   65D25, 65D05, 42C05
% \end{AMS}

\pagestyle{myheadings}

\thispagestyle{plain}

\markboth{M. CARLEY}{Moving least squares via orthogonal polynomials}

\section{Introduction}
\label{sec:intro}

The moving least squares method~\cite[chapter~7]{fasshauer03} is a
technique for interpolation~\cite{levin98} and
differentiation~\cite{marshall-grant-gossler-huyer00,moussa-carley08,%
  belward-turner-ilic08,chenoweth-soria-ooi09,zhang-zhao-liew09,%
  schonauer-adolph01,schonauer-adolph02} on scattered data. The
purpose of this paper is to examine the moving least squares problem
in the framework of orthogonal polynomials, as applied to the
estimation of derivatives.

In applications of the type considered here, the data supplied are the
positions of $N$ points $\bx_{i}$, $i=1,\ldots,N$ and corresponding
values $f_{i}$. At one of these points the derivative is to be
estimated. This is done using an interpolating polynomial
$P(\mathbf{x})$ which minimizes the error:
\begin{align}
  \label{equ:error}
  E &= \sum_{i=1}^{N} w_{i}
  \left(
    P(\bx_{i}) - f_{i}
  \right)^{2},
\end{align}
where $w_{i}$ is a strictly positive weight. The polynomial
$P(\mathbf{x})$ can be computed by direct solution of a least squares
problem and then used to interpolate $f(\bx)$ or to estimate its
derivatives.

There are a number of applications where moving least squares is used
to estimate derivatives of a function specified at discrete
points. One is the estimation of gradients of vorticity in Lagrangian
vortex methods~\cite{marshall-grant-gossler-huyer00,moussa-carley08},
where the gradients are estimated in two or three dimensions 
by fitting a second order polynomial to the components of vorticity
and differentiating the polynomial. It was noted that ``when
computational points become very isolated, due to inadequate spatial
resolution, the condition number of the matrix [used in fitting the
polynomial] becomes very
large''~\cite{marshall-grant-gossler-huyer00}. The solution proposed
for this ill-conditioning was to add additional points to the fit. It
appears that this problem may have been caused by another effect which
has been noted by authors who use moving least squares to solve
partial differential equations on irregular meshes or using mesh-free
techniques. 

In the work of Sch\"{o}nauer and
Adolph~\cite{schonauer-adolph01,schonauer-adolph02}, a finite
difference stencil is developed using polynomials which interpolate
data on points of an unstructured mesh. The points used in the
polynomials are selected by choosing more points than there are
coefficients in the fitting polynomial because ``in $m$ nodes usually
there is not sufficient information for the $m$
coefficients''~\cite{schonauer-adolph01} or, restated, ``there are
linear dependencies on straight lines''~\cite{schonauer-adolph02}. The
number of extra points used in fitting the polynomial was determined
through experience and testing. This raises the issue of the
arrangement of the points used in deriving a polynomial fit.

The issue has been addressed recently by Chenoweth et
al.~\cite{chenoweth-soria-ooi09} who considered the problem of how to
find a least squares fit on points of an unstructured mesh in order to
generate a stencil, while avoiding singularities caused by particular
point configurations, a general form of the problem of ``linear
dependencies''~\cite{schonauer-adolph02}. They state the conditions
under which such singularities will arise and state a criterion
determining when it will not be possible to make a least squares fit
of a given order on a given set of points in two dimensions. This will
happen when selected points are spanned by the same polynomial, for
example, when fitting a second order polynomial to points which lie on
an ellipse in two dimensions. They also give an algorithm for a moving
least squares fit which determines when more points must be added in
order to avoid singularities, and which additional points will be
useful.

Another recent paper employing moving least squares methods for
three-dimensional meshless methods~\cite{zhang-zhao-liew09} proposes
an approach which may help avoid the problem of singularities. The
method is to derive a set of basis functions which are orthogonal with
respect to an inner product defined on the set of points. The use of
orthogonal functions has the advantage of improving the condition
number of the system to be solved to form the least squares fit and,
in this case, allows a smaller number of basis functions to be
used. The authors do not, however, discuss the problem of singular
point configurations other than stating the number of points included
in the fit must be large enough to make the system matrix regular,
which corresponds to the avoidance of singular or ill-conditioned
arrangements of nodes.

Strangely, there does not yet appear to be a published moving least
squares method which explicitly frames the problem in terms of
orthogonal polynomials. The aim of this paper is to present a method
using results from the theory of orthogonal polynomials in multiple
variables~\cite{xu04b} to restate the problem in a manner which
detects singular point configurations and generates a set of
orthogonal polynomials which are unique for the points considered. The
polynomials derived can then be used directly in computing a fit to
the function on the specified data points. The method is quite general
and does not require a knowledge of which configurations give rise to
singularities. In three or more dimensions these configurations are
not easily visualized and, furthermore, a singular value decomposition
becomes increasingly expensive.

\section{Discrete orthogonal polynomials for scattered data}
\label{sec:calc}

The theory of classical orthogonal polynomials of several variables is
well-developed~\cite{xu04b} but that of polynomials orthogonal on
discrete points is not as advanced. A recent paper~\cite{xu04},
however, establishes basic properties of discrete orthogonal
polynomials and gives algorithms for their derivation. In particular,
it establishes the theoretical foundations which allow us to say,
given a set of points, whether orthogonal polynomials of a given order
exist on these points and, if they do, what those polynomials are. In
this section, we will summarize the mathematical tools required to
derive and apply polynomials orthogonal on discrete points. We use the
standard notation in which a polynomial of several variables is
defined as a weighted sum of monomials:
\begin{align}
  \label{equ:polynomial}
  P(\bx) &= \sum_{i=1}^{n} A_{i}\bx^{\ba_{i}},
\end{align}
where $\bx=(x_{1},x_{2},\ldots,x_{d})$, $\bx\in \field{R}^{d}$,
$\ba=(\alpha_{1},\alpha_{2},\ldots,\alpha_{d})$, $\ba\in
\field{N}_{0}^{d}$ and the monomial terms $\bx^{\ba}=\prod_{j=1}^{d}
x_{j}^{\alpha_{j}}$. The degree of $P(\bx)$ is $\max |\ba_{i}|$ where
$|\ba|=\sum_{j=1}^{d}\alpha_{j}$.

% Although the theory of classical orthogonal polynomials of several
% variables is well-developed~\cite{xu04b}, the theory of their discrete
% counterparts is not as advanced~\cite{xu04}. There is, however,
% sufficient theoretical foundation to allow us to derive orthogonal
% polynomials of a given order on a given set of points, if such
% polynomials exist, to establish if such polynomials exist and to apply
% these polynomials to the problems of interpolation and
% differentiation. The theory used is that of Xu~\cite{xu04} who
% summarizes some known results and derives some basic properties of
% orthogonal polynomials on sets of isolated points. The method of this
% paper relies on the application of Xu's method for deriving such
% polynomials.

\subsection{Generation of orthogonal polynomials}
\label{sec:generation}

The first basic tool required is a scheme to generate a set of
polynomials which are orthogonal on a given set of points with respect
to some weight function. This can be done using standard matrix
operations~\cite{galassi-davies-etal05,golub-van-loan96} using the
procedure given by Xu~\cite{xu04}. First, we define the inner product:
\begin{align*}
  \langle f(\bx_{i}), g(\bx_{i}) \rangle
  = \sum_{i=1}^{N} f(\bx_{i})g(\bx_{i})w_{i},
\end{align*}
where $f$ and $g$ are functions evaluated at the data points $\bx_{i}$
and $w_{i}$ is the weight corresponding to $\bx_{i}$, with
$w_{i}>0$.

The first step in generating the orthogonal polynomials is to find a
set of monomial powers $\ba_{j}$ which spans the polynomial space on
$\bx_{i}$. This is done by starting with the monomial $1$ and
systematically adding monomials of increasing degree $\ba_{j}$. 
As each monomial is added, a matrix
\begin{align*}
  X &= 
  \left[
    \begin{matrix}
     \bx_{1}^{\ba_{1}} & \bx_{2}^{\ba_{1}} & \dots &
     \bx_{N}^{\ba_{1}} \\
     \bx_{1}^{\ba_{2}} & \bx_{2}^{\ba_{2}} & \dots &
     \bx_{N}^{\ba_{2}} \\
     \vdots & \vdots & & \vdots \\
     \bx_{1}^{\ba_{n}} & \bx_{2}^{\ba_{n}} & \dots &
     \bx_{N}^{\ba_{n}}
    \end{matrix}
  \right],
\end{align*}
is generated for some initial value $n$. New rows are added to $X$ for
successive values of $\ba_{j}$, taken in lexicographical order at each
$|\ba|$. The rank of $X$ is checked at each step; if it is
rank-deficient, the newly added monomial is rejected. Otherwise, it is
added to the list of $\ba_{j}$ to be included in generating the
polynomials. Rejection of a monomial power will happen because the
point configuration is singular for the combination of monomials which
would result from including the new $\ba_{j}$. Monomials are added
until $X$ is square and of full rank. The output of this procedure is
a list of monomial powers which together span the polynomial space on
the data points.

To generate the orthogonal polynomials from the list of monomials, the
following procedure is used:
\begin{enumerate}
\item generate the symmetric, positive definite matrix $M$, with
  $M_{ij}=\langle \bx^{\ba_{i}}, \bx^{\ba_{j}}\rangle$.
\item perform the decomposition $M=SDS^{T}$, where $D$ is a diagonal
  matrix and $S$ is lower triangular.
\item solve $S^{T}R=D^{-1/2}$ where
  $D^{-1/2}=\text{diag}\{(d_{1}w_{1})^{-1/2},\dots,(d_{N}w_{N})^{-1/2}\}$. This
  can be done using an $LU$ solver with rearrangement of the matrix
  entries.
\item the matrix $R$ now contains the coefficients of the orthogonal
  polynomials. 
\end{enumerate}
In implementing the method, we note that $S^{T}$ can be found directly
by using the algorithm given by Golub and Van
Loan~\cite[page~138]{golub-van-loan96} with the row and column indices
switched.

The orthogonal polynomials $P_{i}$ are now:
\begin{align*}
  P_{i}(\bx) &= \sum_{j=1}^{n} R_{ij}\bx^{\ba_{j}},
\end{align*}
and for later convenience, we scale the coefficients on the inner
products $\langle P_{i}(\bx)P_{i}(\bx)\rangle$ to give an orthonormal
basis.

\subsection{Fitting data on sets of scattered points}
\label{sec:fitting}

Given the set of orthogonal polynomials $P_{i}(\bx)$, generation of a
least-squares fit is trivial. By orthogonality:
\begin{align}
  \label{equ:fit}
  f(\bx) &\approx \sum_{i} c_{i}P_{i}(\bx),
\end{align}
where the constants $c_{i}$ are given by:
\begin{align*}
  c_{i} &= \sum_{j} f_{j}w_{j}P_{i}(\bx_{j}).
\end{align*}
Rearranging to give the interpolant as a weighted sum over the data
points:
\begin{align}
  \label{equ:fit:1}
  f(\bx) &\approx \sum_{j} v_{j}f_{j},\\
  v_{j} &= w_{j}\sum_{i} P_{i}(\bx_{j})P_{i}(\bx).\nonumber
\end{align}
Derivatives of $f(\bx)$ can also be estimated as a weighted sum of the
function values at the points of the distribution, to generate a
differentiation stencil:
\begin{align}
  \label{equ:diff}
  \frac{\partial^{l+m+\cdots}}{\partial^{l}x_{1}\partial^{m}x_{2}\cdots}
  f(\bx) &\approx \sum_{j} v^{(lm\cdots)}_{j}f_{j},\\
  v^{(lm\cdots)}_{j} &= w_{j}\sum_{i} P_{i}(\bx_{j})
  \frac{\partial^{l+m+\cdots}}{\partial^{l}x_{1}\partial^{m}x_{2}\cdots}
  P_{i}(\bx),\nonumber
\end{align}
with the derivatives of $P_{i}(\bx)$ being computed directly from the
coefficients in the matrix $R$ of \S\ref{sec:generation}.

In summary, a derivative of a function $f(\bx)$ given on a set of
points can be estimated at some point $\bx_{0}$ using these steps:
\begin{enumerate}
\item\label{select} select $N$ points in the region of $\bx_{0}$,
  including $\bx_{0}$ itself;
\item\label{ortho} generate a set of orthonormal polynomials for the
  selected points, using the procedure of \S\ref{sec:generation};
\item evaluate the weights $v^{(lm\cdots)}_{j}$ given by
  equation~\ref{equ:diff};
\item calculate the derivative as the weighted sum of the function
  values. 
\end{enumerate}

An important point is that strictly this procedure can only evaluate
linear combinations of derivatives. In an extreme case, where only two
points are used, the available monomials will be $1$ and $x_{1}$ (or
$x_{2}$ depending on the lexicographical ordering used). This allows
linear interpolation of a function $f$ between the two points and
estimation of the derivative on a straight line joining them. This
derivative will be a linear combination of $\partial f/\partial x_{1}$
and $\partial f/\partial x_{2}$, with the precise combination
depending on the orientation of the two points. In practice, this
should not be a serious limitation since the monomials used in the
polynomials are known and it is possible to determine whether there
is a full set available for the determination of all derivatives of a
given order.

\section{Performance}
\label{sec:performance}

To illustrate the operation of the method, the first results presented
are orthogonal polynomials on regular arrangements of points. The
first is a regular $3\times3$ grid in $(-1,1)\times(-1,1)$. Upon
application of the algorithm of \S\ref{sec:generation}, the monomials
which form the final matrix $X$ are $1$, $x_{1}$, $x_{2}$,
$x_{1}^{2}$, $x_{1}x_{2}$, $x_{2}^{2}$, $x_{1}^{2}x_{2}$,
$x_{1}x_{2}^{2}$ and $x_{1}^{2}x_{2}^{2}$ and the resulting orthogonal
polynomials are:
\begin{subequations}
  \label{equ:grid:poly}
  \begin{align}
    P_{0} &= \frac{1}{3},\\
    P_{1} &= \frac{1}{2^{1/2}3^{1/2}}x_{1},\\
    P_{2} &= \frac{1}{2^{1/2}3^{1/2}}x_{2},\\
    P_{3} &= -\frac{2^{1/2}}{3} + \frac{1}{2^{1/2}}x_{1}^{2},\\
    P_{4} &= \frac{1}{2}x_{1}x_{2},\\
    P_{5} &= -\frac{2^{1/2}}{3} + \frac{1}{2^{1/2}}x_{2}^{2},\\
    P_{6} &= -\frac{1}{3^{1/2}}x_{2} + \frac{3^{1/2}}{2}x_{1}^{2}x_{2},\\
    P_{7} &= -\frac{1}{3^{1/2}}x_{1} + \frac{3^{1/2}}{2}x_{1}x_{2}^{2},\\
    P_{8} &= \frac{2}{3} - x_{1}^{2} - x_{2}^{2} + \frac{3}{2}x_{1}^{2}x_{2}.
  \end{align}
\end{subequations}

If the procedure is applied to six points equally spaced on a unit
circle, the resulting polynomials are:
\begin{subequations}
  \label{equ:circle:poly}
  \begin{align}
    P_{0} &= \frac{1}{2^{1/2}3^{1/2}},\\
    P_{1} &= \frac{1}{3^{1/2}}x_{1},\\
    P_{2} &= \frac{1}{3^{1/2}}x_{2},\\
    P_{3} &= -\frac{1}{3^{1/2}} + \frac{2}{3^{1/2}}x_{1}^{2},\\
    P_{4} &= \frac{2}{3^{1/2}}x_{1}x_{2},\\
    P_{5} &= -\frac{3^{1/2}}{2^{1/2}}x_{1} +
    2\frac{2^{1/2}}{3^{1/2}}x_{1}^{3}. 
  \end{align}
\end{subequations}

A number of general issues are illustrated by these examples. The
first is the obvious one that there are no more polynomials than there
are points. This means that although the polynomials are notionally up
to third order in both cases, in practice neither group of functions
has a complete set of monomials capable of spanning all polynomials up
to cubic. Secondly, if the orthogonal polynomials can be generated,
there is no benefit in adding extra points once a complete set of
functions is available: the result of adding more points is to yield
an incomplete set of higher order polynomials. In applications, it may
well be better to have a lower order, but complete, system to fit the
functions on the points.

% \begin{itemize}
% \item high order on too few points causes problems due to dropped terms;
% \item low order does not benefit from having more points.
% \end{itemize}

\subsection{Random point distributions on the unit disc or ball}
\label{sec:random}

The first set of results presented are average data for a large number
of tests conducted with varying order and length scale. Following the
example of Belward et al.~\cite{belward-turner-ilic08}, the accuracy
and robustness of the computational scheme are tested by estimating
the derivatives of a prescribed function using a set of points
randomly distributed over the unit disc or ball. The functions used
are:
\begin{subequations}
  \label{equ:test}
  \begin{align}
    \label{equ:poly}
    f_{1}(\mathbf{x}) &= R^{4},\\
    f_{2}(\mathbf{x}) &= \E^{-R^{2}},\\
    f_{3}(\mathbf{x}) &= x_{1}\E^{-R^{2}},\\
%    f_{4}(\mathbf{x}) &= \frac{\sin \gamma R}{R},\\
    R^{2} &= \sum x_{j}^{2}.\nonumber
  \end{align}
\end{subequations}
The functions have been chosen to give a function which can be fitted
exactly by a polynomial ($f_{1}$), a Gaussian of the type found in
various applications such as vortex dynamics ($f_{2}$) and a Gaussian
weighted to give an asymmetry with a consequent non-zero derivative at
the evaluation point ($f_{3}$). The evaluation point was fixed at
$\mathbf{0}$ and $N=8,16,32,64,128$ random points were distributed
uniformly in angle and radius over the unit disc or ball. Unit weights
$w_{i}\equiv1$ were used. Given the values of $f(\mathbf{x}_{i})$, the
algorithm of~\S\ref{sec:fitting} was used to estimate first and second
derivatives of the function at the evaluation point. To examine the
convergence rate, the procedure was repeated by using the values
$f(\sigma x_{i})$, where $\sigma$, $0<\sigma\leq 1$, is a scaling
factor which has the effect of contracting the point distribution
around the evaluation point. It is assumed that in applications, a
point distribution will be scaled to a normalized radius and the
result of the function evaluation rescaled afterwards, a procedure
which is modelled using the scaling factor $\sigma$. Tests were
repeated on~32 different random point distributions and mean absolute
errors estimated. The presented results are mean absolute error, mean
number of rejected monomials and convergence rate with $\sigma$, for
different values of $N$ and different functions $f(\mathbf{x}$ in two
and three dimensions.
%A tolerance of $10^{-3}$ was used in the test of the 

Two sample sets of results are shown graphically to illustrate the
performance of the method, with the relevant performance parameters
for all tests summarized in tabular form. 

\begin{figure}
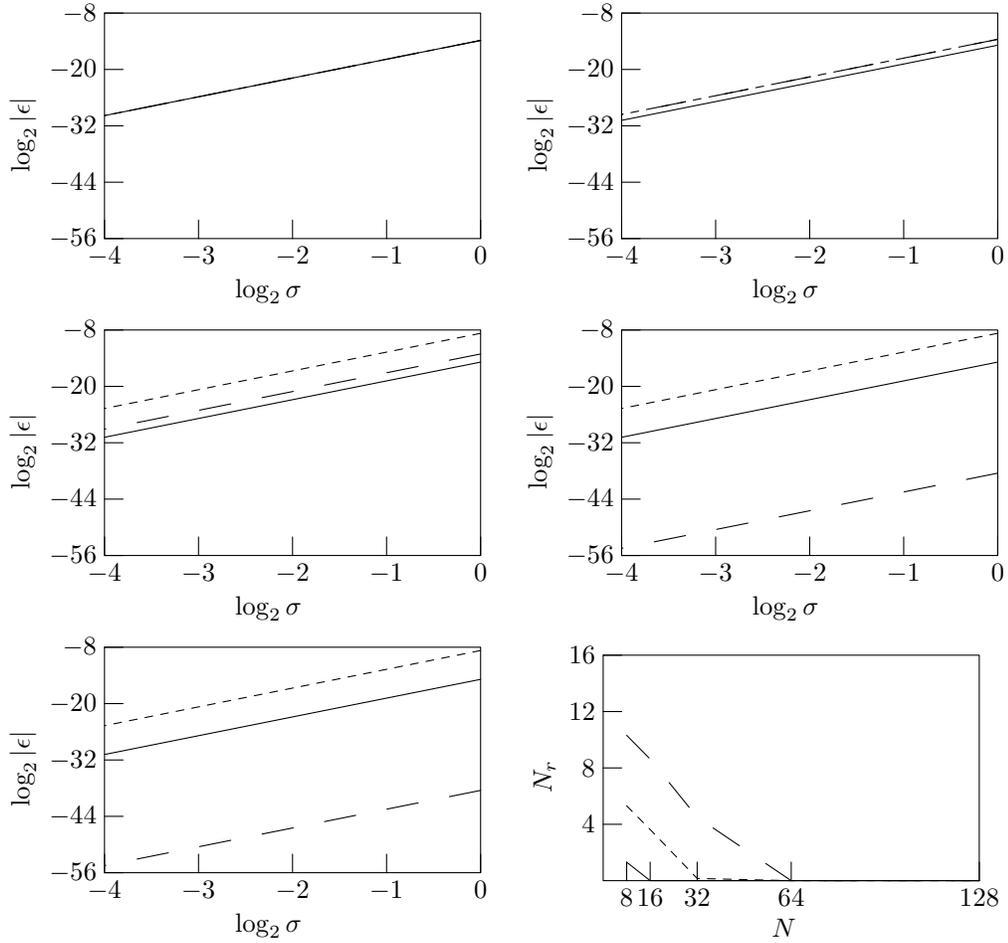

  \centering
  \begin{tabular}{cc}
    \includegraphics{siam07a-figs.1} &
    \includegraphics{siam07a-figs.2} \\

    \includegraphics{siam07a-figs.3} &
    \includegraphics{siam07a-figs.4} \\
    
    \includegraphics{siam07a-figs.5} &
    \includegraphics{siam07a-figs.6} 
    
  \end{tabular}

  \caption{Evaluation of $\partial f_{1}/\partial x_{1}$ in two
    dimensions: mean absolute error against scale factor $\sigma$ for, reading
    left to right, $N=8,16,32,64,128$ points; final plot mean number of
    rejected monomials $N_{r}$ against number of points $N$. Second
    order fit shown solid, third order dashed, fourth order long dashed.}
  \label{fig:err:f1:2d}
\end{figure}

Figure~\ref{fig:err:f1:2d} shows the performance data for evaluation
of $\partial f_{1}/\partial x_{1}$ in two dimensions. This function
can be fitted exactly by a polynomial of sufficiently high order as is
clear from the results. The first plot of mean absolute error against
$\sigma$, shows that all three orders give identical results. This is
because, with only eight points available in the fit, the three
systems of functions are identical. Increasing the number of points to
sixteen, the second order fit is slightly better than the other two
which are themselves identical. The final plot, showing the mean
number of monomials rejected for each $N$ explains why. With $N=8$,
all three fits are underspecified while at $N=16$, the second order
fit has a full set of six monomials but the third and fourth order
fits are still short of useable terms. As $N$ increases to~64, the
third and fourth order fits gain a complete set of monomials (ten and
fifteen, respectively) and the full accuracy of the exact fourth order
fit becomes clear. The third order fit's error, however, is larger
than that of the second order, probably because the third order
monomials, which do not fit the symmetric function $f_{1}$, introduce
spurious terms in the fit. The results for $N=64$ and $N=128$ are
identical, because at this stage a full set of monomials is available
for all of the orders considered and adding additional points gives no
extra benefit. In all of the cases considered, the error reduces with
$\sigma$ at the same rate, though with quite different error
magnitudes, as will be seen in the tabular data presented later.

\begin{figure}
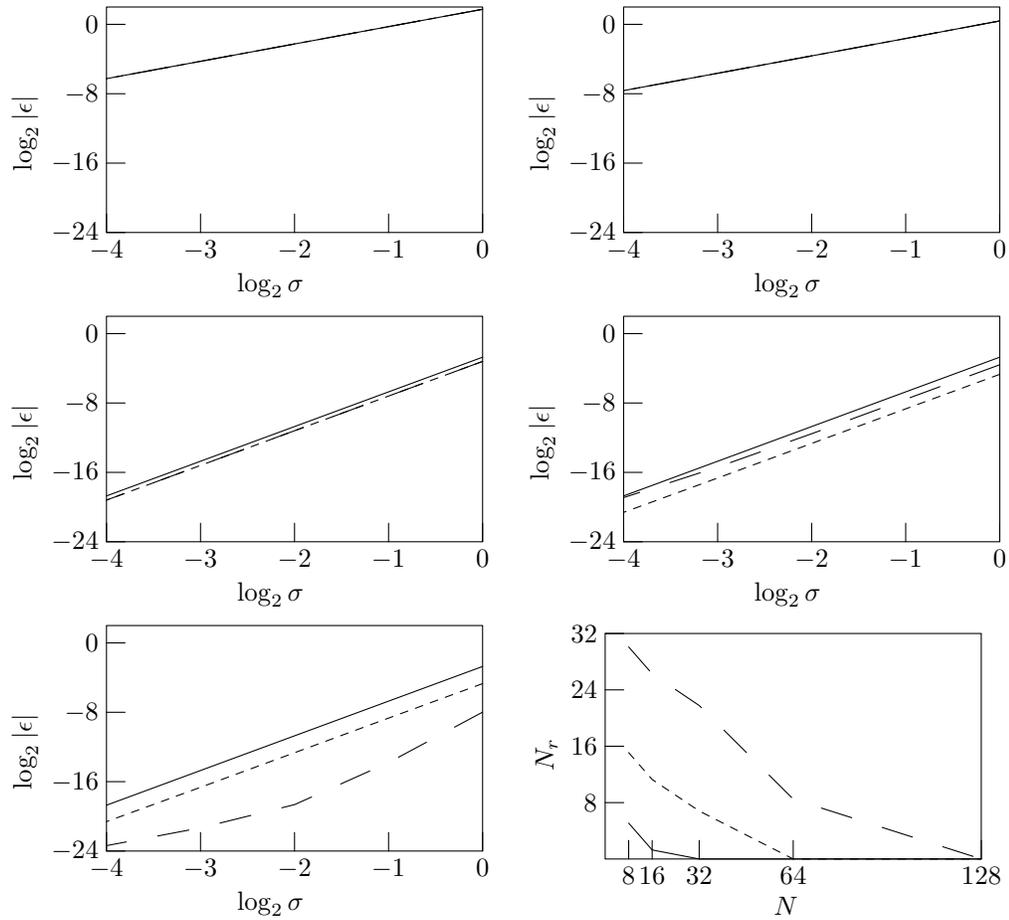

  \centering
  \begin{tabular}{cc}
    \includegraphics{siam07a-figs.11} &
    \includegraphics{siam07a-figs.12} \\

    \includegraphics{siam07a-figs.13} &
    \includegraphics{siam07a-figs.14} \\
    
    \includegraphics{siam07a-figs.15} &
    \includegraphics{siam07a-figs.16} 
    
  \end{tabular}

  \caption{Evaluation of $\partial^{2} f_{2}/\partial x_{1}^{2}$ in
    three dimensions: mean absolute error $\bar{|\epsilon|}$ against
    scale factor $\sigma$ for, reading left to right,
    $N=8,16,32,64,128$ points; final plot mean number of rejected
    monomials $N_{r}$ against number of points $N$. Second order fit
    shown solid, third order dashed, fourth order long dashed.}
  \label{fig:err:f2:3d}
\end{figure}

Figure~\ref{fig:err:f2:3d} shows the performance of evaluation of
$\partial f_{2}/\partial x_{1}^{2}$ in three dimensions. For
reference, there are~10,~20 and~35 monomials in a fully-specified
polynomial of order two, three and four respectively. As in the
two-dimensional case, the first two plots show identical behavior of
the error for all three fits, due to the number of points being
insufficient to generate a fully specified polynomial. The third plot,
$N=32$, shows the second order fit being slightly better than the
other two, which are identical, since this is now fully specified. As
the number of points is increased the three fits begin to separate,
although none has a clear advantage over the others until $N=128$,
where the fourth order fit has a full set of monomials available, as
shown in the final plot. The convergence rate, however, reduces at
small $\sigma$ which may be due to floating point errors. In the
fourth order fit the monomials are $O(x^{4})$ and the resulting inner
products $O(x^{8})$. When $\sigma=2^{-4}$, the maximum value of
$x^{8}$ is $2^{-32}\approx2\times10^{-10}$, for points furthest from
the evaluation position: the corresponding term for those points
nearest the evaluation position will be much smaller, comparable to
the floating point precision of the computer. 

\begin{table}
  \caption{Convergence rates for first derivatives in two-dimensional
    problems}
  \label{tab:conv:f1:2d}
  \centering
  \begin{tabular}{lcrrrrrrr}
    & & 
    \multicolumn{1}{c}{8} &
    \multicolumn{1}{c}{16} &
    \multicolumn{1}{c}{32} &
    \multicolumn{1}{c}{64} &
    \multicolumn{1}{c}{128} &
    \multicolumn{1}{c}{$\epsilon_{\min}^{(128)}$} &
    \multicolumn{1}{c}{$\epsilon_{\max}^{(128)}$} 
    \\
    \hline
    $f_{1}$  & 2 & 4.00 & 4.00 & 4.00 & 4.00 & 4.00& 5.21$\times10^{-10}$ & 3.41$\times10^{-5}$\\
 & 3 & 4.00 & 4.00 & 4.00 & 4.00 & 4.00& 3.66$\times10^{-8}$ & 2.40$\times10^{-3}$\\
 & 4 & 4.00 & 4.00 & 4.00 & 4.00 & 4.00& 3.98$\times10^{-17}$ & 2.61$\times10^{-12}$\\
$f_{2}$  & 2 & 2.00 & 2.03 & 4.00 & 4.00 & 4.00& 1.09$\times10^{-9}$ & 7.18$\times10^{-5}$\\
 & 3 & 2.00 & 2.05 & 4.00 & 4.00 & 4.00& 3.03$\times10^{-9}$ & 1.98$\times10^{-4}$\\
 & 4 & 2.00 & 2.05 & 4.00 & 5.42 & 5.42& 4.99$\times10^{-13}$ & 1.24$\times10^{-6}$\\
$f_{3}$  & 2 & 3.00 & 3.00 & 3.00 & 3.00 & 3.00& 2.18$\times10^{-7}$ & 8.93$\times10^{-4}$\\
 & 3 & 3.00 & 3.00 & 3.02 & 4.80 & 4.80& 1.96$\times10^{-11}$ & 1.05$\times10^{-5}$\\
 & 4 & 3.00 & 3.00 & 3.02 & 4.35 & 4.35& 4.99$\times10^{-11}$ & 6.22$\times10^{-6}$\\

  \end{tabular}
\end{table}

Numerical data summarizing the results of all of the tests carried out
are presented in Tables~\ref{tab:conv:f1:2d}--\ref{tab:conv:f2:3d}. In
each table, the convergence rate of the fits is presented for each of
the three test functions at each value of $N$, the size of the point
set. In addition, to compare the errors proper, the final two columns
give the maximum and minimum mean errors for fits performed using~128
points. Convergence rates $r$ were found by a least squares fit
$|\epsilon|=\epsilon_{0}\sigma^{r}$.

Table~\ref{tab:conv:f1:2d} shows the performance data for evaluation
of $\partial f_{i}/\partial x_{1}$. The results are much as might be
expected, with a small minimum error in each case, especially for
$f_{1}$, the polynomial and with smooth convergence for most fits on
most functions. The exceptions are the fourth order fit to $f_{2}$ and
$f_{3}$ and the third order fit to $f_{3}$. In these cases, at large
point number $N\geq64$, the errors are small, as would be expected,
but the convergence is not as smooth as expected. For $f_{3}$ the mean
convergence rate of the fourth order fit is also less than that of the
third order, although the absolute errors are comparable. This is
probably due to a combination of the floating point issue mentioned
earlier and the inability of the fourth order monomials to capture the
behaviour of the function near $\mathbf{0}$. 

\begin{table}
  \caption{Convergence rates for second derivatives in two-dimensional
    problems}
  \label{tab:conv:f2:2d}
  \centering
  \begin{tabular}{lcrrrrrrr}
    & & 
    \multicolumn{1}{c}{8} &
    \multicolumn{1}{c}{16} &
    \multicolumn{1}{c}{32} &
    \multicolumn{1}{c}{64} &
    \multicolumn{1}{c}{128} &
    \multicolumn{1}{c}{$\epsilon_{\min}^{(128)}$} &
    \multicolumn{1}{c}{$\epsilon_{\max}^{(128)}$}
    \\
    \hline
    $f_{1}$  & 2 & 4.00 & 4.00 & 4.00 & 4.00 & 4.00& $7.39\times10^{-8}$ & $4.84\times10^{-3}$\\
 & 3 & 4.00 & 4.00 & 4.00 & 4.00 & 4.00& $1.34\times10^{-7}$ & $8.84\times10^{-3}$\\
 & 4 & 4.00 & 4.00 & 4.00 & 4.00 & 4.00& $2.65\times10^{-15}$ & $1.74\times10^{-10}$\\
$f_{2}$  & 2 & 2.00 & 2.02 & 4.00 & 4.00 & 4.00& $6.25\times10^{-8}$ & $4.09\times10^{-3}$\\
 & 3 & 2.00 & 2.03 & 4.00 & 4.00 & 4.00& $4.10\times10^{-6}$ & $2.67\times10^{-1}$\\
 & 4 & 2.00 & 2.03 & 4.00 & 5.01 & 5.01& $7.49\times10^{-11}$ & $6.02\times10^{-5}$\\
$f_{3}$  & 2 & 3.00 & 3.00 & 3.00 & 3.00 & 3.00& $3.15\times10^{-5}$ & $1.29\times10^{-1}$\\
 & 3 & 3.00 & 3.00 & 3.00 & 4.88 & 4.88& $1.37\times10^{-10}$ & $9.64\times10^{-5}$\\
 & 4 & 3.00 & 3.00 & 3.01 & 4.75 & 4.75& $3.29\times10^{-9}$ & $1.52\times10^{-3}$\\

  \end{tabular}
\end{table}

Table~\ref{tab:conv:f2:2d} shows the equivalent results for evaluation
of $\partial^{2} f_{i}\partial x_{1}^{2}$. The results show the same
trends as in Table~\ref{tab:conv:f1:2d}, with the fourth order fit
giving very small minimum errors for all three functions but with the
third order fit being slightly superior for $f_{3}$. 

\begin{table}
  \caption{Convergence rates for first derivatives in three-dimensional
    problems}
  \label{tab:conv:f1:3d}
  \centering
  \begin{tabular}{lcrrrrrrr}
    & & 
    \multicolumn{1}{c}{8} &
    \multicolumn{1}{c}{16} &
    \multicolumn{1}{c}{32} &
    \multicolumn{1}{c}{64} &
    \multicolumn{1}{c}{128} &
    \multicolumn{1}{c}{$\epsilon_{\min}^{(128)}$} &
    \multicolumn{1}{c}{$\epsilon_{\max}^{(128)}$}
    \\
    \hline
    $f_{1}$  & 2 & 4.00 & 4.00 & 4.00 & 4.00 & 4.00& $3.14\times10^{-7}$ & $2.05\times10^{-2}$\\
 & 3 & 4.00 & 4.00 & 4.00 & 4.00 & 4.00& $1.73\times10^{-8}$ & $1.13\times10^{-3}$\\
 & 4 & 4.00 & 4.00 & 4.00 & 4.00 & 4.00& $4.55\times10^{-14}$ & $2.98\times10^{-9}$\\
$f_{2}$  & 2 & 2.00 & 2.00 & 4.00 & 4.00 & 4.00& $4.73\times10^{-9}$ & $3.09\times10^{-4}$\\
 & 3 & 2.00 & 2.00 & 3.97 & 4.00 & 4.00& $1.18\times10^{-6}$ & $7.79\times10^{-2}$\\
 & 4 & 2.00 & 2.00 & 3.97 & 3.99 & 3.61& $9.13\times10^{-10}$ & $2.31\times10^{-5}$\\
$f_{3}$  & 2 & 3.00 & 3.00 & 3.00 & 3.00 & 3.00& $5.12\times10^{-7}$ & $2.09\times10^{-3}$\\
 & 3 & 3.00 & 3.00 & 3.00 & 4.75 & 4.75& $1.50\times10^{-10}$ & $6.89\times10^{-5}$\\
 & 4 & 3.00 & 3.00 & 3.00 & 2.01 & 4.39& $7.88\times10^{-9}$ & $1.16\times10^{-3}$\\

  \end{tabular}
\end{table}

\begin{table}
  \caption{Convergence rates for second derivatives in three-dimensional
    problems}
  \label{tab:conv:f2:3d}
  \centering
  \begin{tabular}{lcrrrrrrr}
    & & 
    \multicolumn{1}{c}{8} &
    \multicolumn{1}{c}{16} &
    \multicolumn{1}{c}{32} &
    \multicolumn{1}{c}{64} &
    \multicolumn{1}{c}{128} &
    \multicolumn{1}{c}{$\epsilon_{\min}^{(128)}$} &
    \multicolumn{1}{c}{$\epsilon_{\max}^{(128)}$}  
    \\
    \hline
    $f_{1}$  & 2 & 4.00 & 4.00 & 4.00 & 4.00 & 4.00& $1.66\times10^{-6}$ & $1.09\times10^{-1}$\\
 & 3 & 4.00 & 4.00 & 4.00 & 4.00 & 4.00& $1.35\times10^{-6}$ & $8.91\times10^{-2}$\\
 & 4 & 4.00 & 4.00 & 4.00 & 4.00 & 4.00& $1.11\times10^{-11}$ & $7.32\times10^{-7}$\\
$f_{2}$  & 2 & 2.00 & 2.01 & 4.00 & 4.00 & 4.00& $2.32\times10^{-6}$ & $1.51\times10^{-1}$\\
 & 3 & 2.00 & 2.01 & 4.00 & 3.98 & 3.98& $6.25\times10^{-7}$ & $3.85\times10^{-2}$\\
 & 4 & 2.00 & 2.01 & 4.00 & 3.85 & 3.83& $9.17\times10^{-8}$ & $3.95\times10^{-3}$\\
$f_{3}$  & 2 & 3.00 & 3.00 & 3.00 & 3.00 & 3.00& $7.14\times10^{-5}$ & $2.91\times10^{-1}$\\
 & 3 & 3.00 & 3.00 & 3.00 & 5.00 & 5.00& $4.23\times10^{-9}$ & $4.44\times10^{-3}$\\
 & 4 & 3.00 & 3.00 & 3.00 & 4.13 & 3.50& $3.36\times10^{-6}$ & $4.56\times10^{-2}$\\

  \end{tabular}
\end{table}

Tables~\ref{tab:conv:f1:3d} and~\ref{tab:conv:f2:3d} give data for
first and second derivative evaluation in three dimensions. The
general trends are similar to those in two dimensions but there some
differences worth noting. In evaluating the derivatives of $f_{3}$,
the fourth order fit, as in the two-dimensional case, does not perform
as well as the third order, although the error is still small. The
difference is in the convergence rate as $N$ increases from~32. At
$N=64$, the convergence rate drops from~3 to~2 before increasing again
to~4.39, in contrast to the behavior of the third order fit. This is
probably because at $N=32$, the third and fourth order fits are both
underspecified (see the final plot of Figure~\ref{fig:err:f2:3d}) but
the third order fit gains a full set of monomials at $N=64$. The
fourth order fit has a full set of third order monomials but has
rejected, on average, eight monomials, leaving only seven fourth order
terms in the polynomials. This leads to error behavior which is worse
than that from a fully specified third order polynomial. As $N$
increases again, to~128, the fourth order polynomial is fully-defined
and its convergence rate recovers, with the caveat that it may now be
affected by floating point errors.

\begin{figure}
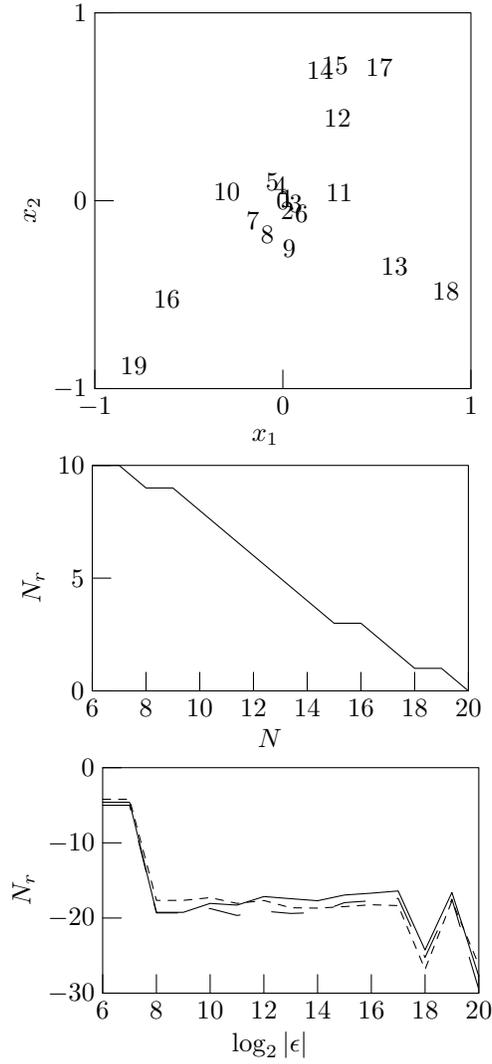

  \centering
  \includegraphics{siam07a-figs.17} \\
  \includegraphics{siam07a-figs.18} \\
  \includegraphics{siam07a-figs.19}
  \caption{Single point configuration test: top: point configuration;
    middle: number of rejected monomials against number of points
    used; bottom: error in $\partial^{2}f/\partial x_{1}^{2}$ (solid),
    $\partial^{2}f/\partial x_{1}\partial x_{2}$ (dashed) and
    $\partial^{2}f/\partial x_{2}^{2}$ (long dashed).}
  \label{fig:detail}
\end{figure}

To examine the behavior of the method in more detail, we look at the
results of a test on a fixed random point configuration. The task is
to estimate the second derivatives of $f_{2}=\exp(-R^{2})$ in two
dimensions with $\sigma=1/2^{3}$ using a nominally fourth order
fit. Twenty randomly distributed points, sorted by distance from the
evaluation point $\mathbf{x}=\mathbf{0}$, are used and a fit is
generated using $N=6,7,\ldots,20$ of these points in turn. The results
are shown in Figure~\ref{fig:detail}: the first plot shows the point
configuration with points numbered by distance from the evaluation
point, the second plot shows the number of monomials rejected at each
$N$ and the third shows the error in the estimate of each of the three
second-order derivatives.

The error behavior demonstrates some of the detailed features of using
the estimation scheme. For a fully specified fourth order polynomial,
fifteen monomials are required. It is only at $N=20$ that these all
become available, with sufficient points being used to avoid singular
configurations. For $N=6,7$, the error is quite large, due to the fit
being a set of defective second order polynomials---ten monomials are
rejected and only five terms are available for a notionally second
order fit. At $N=8$, one more monomial becomes available and the error
drops immediately since this the method is now a fully specified
second order fit to a set of points close to the evaluation point. The
error remains constant up to $N=18$, even though the number of points
is increasing, since the additional terms are third order and do not
contribute to a fit on the symmetric function in question. There is a
drop in the error at $N=18$, followed by an increase and another
reduction. This can be explained by looking at which monomials have
been included or rejected in the fits.

\begin{table}
  \centering
  \caption{Fourth order monomials included in polynomial fits on point
    set of Figure~\ref{fig:detail} and resulting error in
    $\partial^{2}f/\partial x_{1}\partial x_{2}$}
  \label{tab:detail}
  \begin{tabular}{lcccccc}
    $N$ & \multicolumn{5}{c}{Fourth order monomials included} &
    $\log_{2}|\epsilon|$\\
    \hline
    17 & 
    $x_{1}^{4}$ &
    $x_{1}^{3}x_{2}^{1}$ &
    &
    &
    $x_{2}^{4}$ 
    &
    -18.3
    \\
    18 & 
    $x_{1}^{4}$ &
    $x_{1}^{3}x_{2}^{1}$ &
    $x_{1}^{2}x_{2}^{2}$ &
    &
    $x_{2}^{4}$ 
    &
    -26.8
    \\
    19 & 
    $x_{1}^{4}$ &
    $x_{1}^{3}x_{2}^{1}$ &
    &
    $x_{1}^{1}x_{2}^{3}$ &
    $x_{2}^{4}$ 
    &
    -17.7
    \\
    20 & 
    $x_{1}^{4}$ &
    $x_{1}^{3}x_{2}^{1}$ &
    $x_{1}^{2}x_{2}^{2}$ &
    $x_{1}^{1}x_{2}^{3}$ &
    $x_{2}^{4}$ &
    -26.2
  \end{tabular}
\end{table}

Table~\ref{tab:detail} shows the fourth order terms which are included
in the polynomial fits whose error behavior is shown in
Figure~\ref{fig:detail}, for $17\leq N\leq 20$, with the final column
showing the error in a second order derivative. Each of the fits has a
full set of lower order monomials and, in principle increasing $N$
allows more fourth order terms to be added. In practice, as can be
seen, on this point set, at $N=17$, two fourth order terms have been
rejected and the accuracy suffers. At $N=18$, the term
$x_{1}^{2}x_{2}^{2}$ is added, and the accuracy improves markedly:
this monomial is symmetric and is useful in a fit to the symmetric
function being differentiated. At $N=19$, the monomial
$x_{1}^{2}x_{2}^{2}$ is rejected and $x_{1}^{1}x_{2}^{3}$ is
included. This term increases the error as in seen in the final column
of Table~\ref{tab:detail} and in Figure~\ref{fig:detail}. A full set
of monomials is only available for this point configuration when
$N=20$, resulting in the full expected accuracy. 

This behavior is slightly unexpected given that the orthogonal
polynomials derived for any value of $N$ span the polynomial function
space on those points. This raises the issue of which basis functions
should be used in applications. Chenoweth et
al.~\cite{chenoweth-soria-ooi09} discuss the problem of singular point
configurations in the context of the singular value decomposition of a
matrix containing the monomials evaluated at each point. As these
authors note, a singular value decomposition shows which basis
functions span the null space of polynomials on the points, allowing
the detection of singular point configurations. The opposite fact is
also true: the singular value decomposition yields a set of basis
functions which span the function space on the points and, indeed,
will indicate which of these basis functions are best determined. The
problem, as we see above, is that even when a full set of
well-determined basis functions is available, it is not guaranteed
that they form a suitable basis for the evaluation of derivatives,
unless some extra measures are taken, as in Chenoweth et al's work
where derivatives are included in the general form of the function to
be fitted~\cite{chenoweth-soria-ooi09}. 

Given that in many applications, it will not be known in advance which
terms will be most useful in fitting a function, it is recommended
that only fully specified polynomials be employed, with the order
depending on the accuracy required and the point density available.

\section{Conclusions}
\label{sec:conclusions}

A method for moving least squares interpolation and differentiation
using orthogonal polynomials has been presented and tested on random
point distributions. The method makes use of the theory of discrete
orthogonal polynomials in multiple variables and deals with the
problems caused by singular point configurations by adjusting the
terms of the polynomial. It is concluded that the method is robust and
capable of detecting and compensating for singular configurations. In
applications, it is recommended that the highest order of polynomial
for which a full set of monomials is available be used in computing
derivatives. 

\section*{Acknowledgements}

The author thanks the anonymous referees who read the original
submission with great care, making many useful comments on the method
of the paper and on its presentation.

\bibliography{abbrev,misc,vortex,turbulence,jets,maths,%
  propnoise,identification}

\end{document}